\documentclass[11pt]{amsart}
\usepackage[matrix, arrow]{xy}
\theoremstyle{plain}

\newtheorem{thm}{Theorem}
\newtheorem{theorem}[thm]{Theorem}

\newtheorem{corollary}[thm]{Corollary}
\newtheorem{lemma}[thm]{Lemma}

\newtheorem*{conjecture*}{Conjecture}

\newtheorem*{question*}{Question}

\theoremstyle{definition}

\newtheorem*{rem*}{Remark}
\newtheorem{remark}[thm]{Remark}
\newtheorem*{remark*}{Remark}

\newtheorem*{remarks*}{Remarks}

\newtheorem*{example*}{Example}

\newtheorem*{examples*}{Examples}

\newtheorem*{notation*}{Notation}

\newtheorem*{bibliographical-note}{Bibliographical note}
\newcommand{\zeroindent}{\parindent0cm \parskip1ex}

{\begin{list}{(\alph{enumi}) }{\usecounter{enumi}

\leftmargin0cm \labelsep0cm \rightmargin0cm 
\setlength{\labelwidth}{\fill}}}
{\end{list}}

{

\begin{enumerate}}{\end{enumerate}}

{

\begin{enumerate}}{\end{enumerate}}

{

\begin{enumerate} \parsep0cm 
\itemsep0cm \parskip0cm}{\end{enumerate}}

{

\begin{enumerate}}{\end{enumerate}}

{

\begin{enumerate}}{\end{enumerate}}

\newcommand{\R}{\mathbb{R}}
\newcommand{\Z}{\mathbb{Z}}

\newcommand{\id}{\mathrm{id}}

\newcommand{\iso}{\cong}           %isomorphism sign
\newcommand{\smooth}{C^\infty}

\newcommand{\half}{{\textstyle\frac{1}{2}}}

\renewcommand{\o}{\omega}

\newcommand{\Diff}{\mathrm{Diff}}

\CompileMatrices
\zeroindent

\renewcommand{\SS}{{\mathcal S}}
\newcommand{\DD}{{\mathcal D}}
\newcommand{\GG}{{\mathcal G}}
\newcommand{\F}{{\mathcal F}}
\newcommand{\cotangent}{T^{*\!}S^2}
\newcommand{\Map}{\mathrm{Map}}
\newcommand{\Dbar}{\overline{\Delta}}
\newcommand{\point}{point}
\newcommand{\J}{{\mathcal J}}

\title{Symplectic automorphisms of $\cotangent$}
\author{Paul Seidel}
\date{19/3/1998}
\thanks{Research supported by NSF grant DMS 9304580.}
\address{Institute for Advanced Study, Princeton}
\email{pseidel@math.ias.edu}

\begin{document}

\begin{abstract}
Let $T \subset \cotangent$ be the bundle of unit discs, which is a
compact symplectic manifold with boundary. Let $\SS$ be the group of symplectic
automorphisms of $T$ which are trivial near $\partial T$. In \cite{seidel98b}
it was proved that $\pi_0(\SS)$ contains an element of infinite order. This element
is given by the `generalized Dehn twist' $\tau$ which is the monodromy map 
of a quadratic singularity \cite{arnold95} (see also \cite{seidel97}). The fact
that $[\tau]$ has infinite order is interesting because the class of $\tau$ in
$\pi_0(\DD)$, where $\DD$ is the corresponding group of diffeomorphisms, has
order two. This note contains a short
direct proof of the fact that $\SS$ is weakly homotopy equivalent to the 
discrete group $\Z$, with $\tau$ as a generator. The proof proceeds by
compactifying $T$ to $S^2 \times S^2$. The topology of
the symplectic automorphism group of $S^2 \times S^2$ was determined by
Gromov \cite{gromov85}, and we use a variant of his argument.
\end{abstract}
\maketitle

Let $T \subset \cotangent$ be the space of cotangent vectors of length
$\leq 1$, and $\eta$ its standard symplectic structure. Let $\DD$ be the group 
of diffeomorphisms $\phi: T \longrightarrow T$ such that $\phi = \id$ in a
neighbourhood of $\partial T$, with the $\smooth$-topology, and $\SS$ the
subgroup of those $\phi$ which are symplectic (that is, they preserve
$\eta$).

\begin{theorem} \label{th:ss}
$\pi_0(\SS) \iso \Z$, and $\pi_i(\SS) = 1$ for all $i>0$.
\end{theorem}

The components of $\SS$ cannot all be distinguished in $\DD$:

\begin{corollary} \label{th:ss-dd}
The image of the map $\pi_0(\SS) \longrightarrow \pi_0(\DD)$ is isomorphic
to $\Z/2$.
\end{corollary}

We will prove both results by compactifying $(T,\eta)$ to $S^2 \times S^2$ with its
standard symplectic structure $\o$ (the one for
which both spheres have equal areas). More precisely, one can identify (up to
rescaling the symplectic forms) the interior of $(T,\eta)$ with the complement
of the diagonal $\Delta$ in $(S^2 \times S^2,\o)$. This has the following
consequence: let $\DD_2$ be the group of diffeomorphisms $\psi$ of $S^2 \times S^2$
such that $\psi|\Delta = \id$ and which act trivially on the normal bundle
$\nu\Delta$. Then $\DD_2$ is w.h.e. (weakly homotopy equivalent) to $\DD$.
Similarly, the subgroup $\SS_2 \subset \DD_2$ of those $\psi$ which preserve
$\o$ is w.h.e. to $\SS$.

Consider the larger group $\DD_1 \supset \DD_2$ of oriented diffeomorphisms
of $S^2 \times S^2$ which are equal to the identity on $\Delta$, and the
subgroup $\SS_1 \subset \DD_1$ of those which are symplectic. Let $\GG$ be
the group of gauge transformations of $\nu\Delta$ as an oriented
vector bundle, and $\GG^s \subset \GG$ the subgroup of symplectic gauge
transformations. These topological groups form a commutative diagram with
exact rows:
\begin{equation} \label{eq:diagram}
\xymatrix{
1 \ar[r] &
\SS_2 \ar[r] \ar[d] &
\SS_1 \ar[r] \ar[d] &
\GG^s \ar[r] \ar[d] &
1 \\
1 \ar[r] &
\DD_2 \ar[r] &
\DD_1 \ar[r] &
\GG \ar[r] &
1
}
\end{equation}

\begin{lemma} \label{th:gg}
Let $e_x: \GG^s \longrightarrow SL_2(\R)$ be the evaluation at some
point $x \in \Delta$. $e_x$ is a w.h.e. \end{lemma}

\proof Let $\GG^u \subset \GG^s$ be the subgroup of unitary gauge
transformations. $\GG^u$ is a deformation retract of $\GG^s$, and
$\GG^u_0 = \GG^u \cap \ker(e_x)$ is a deformation retract of $\ker(e_x)$.
Now $\GG^u_0$ is homeomorphic to the space $\Map_b(S^2,S^1)$ of
smooth based maps $S^2 \longrightarrow S^1$. Therefore
$\pi_i(\GG^u) \iso \pi_0(\Map_b(S^{2+i},S^1)) \iso H^1(S^{2+i}) = 0$
for all $i \geq 0$. This proves that $\GG^u_0$, and hence also
$\ker(e_x)$, is weakly contractible. \qed

\begin{lemma} \label{th:ss-one}
Let $\iota \in \SS_1$ be the involution of $S^2 \times S^2$ which
interchanges the two factors. Then the discrete subgroup $\{\id,\iota\} \subset
\SS_1$ is a deformation retract. \end{lemma}

We postpone the proof of this Lemma to the end of the note. Note that the two
components of $\SS_1$ are distinguished by the action on homology.

\proof[Proof of Theorem \ref{th:ss}] We use the long exact sequence of homotopy groups
induced by the top row of \eqref{eq:diagram}. Lemma \ref{th:gg} and Lemma \ref{th:ss-one}
imply that $\pi_i(\SS_2) = 1$ for $i>1$. For $i = 0$ one obtains a short exact sequence
\begin{equation} \label{eq:exact-sequence}
1 \longrightarrow \pi_1(\GG^s) \iso \Z
\stackrel{\partial}{\longrightarrow} \pi_0(\SS_2)
\longrightarrow \pi_0(\SS_1) \iso \Z/2 \longrightarrow 1. 
\end{equation}

Let $\Dbar = \{ (x,-x) \; : \; x \in S^2 \} \subset S^2 \times S^2$ be the antidiagonal.
Its complement carries the Hamiltonian $S^1$-action 
$\rho(t)(x,y) = (R^t_{x+y}(x),R^t_{x+y}(y))$, where $R^t_{\xi} \in SO(3)$ denotes
the rotation with axis $\xi/|\xi|$ and angle $t$. The moment map of $\rho$ is
$\mu(x,y) = -|x+y|$. Choose an $r \in \smooth(\R,\R)$ with $r(t) = -\pi$ for
$t \leq 1/2$ and $r(t) = 0$ for $t \geq 1$. Define $\tau \in \SS_2$ by
\[
\tau(x,y) =
\begin{cases}
(y,x) & |x+y| \leq \half,\\
\rho(r(|x+y|))(x,y) &\text{otherwise.}
\end{cases}
\]
Since $\tau|\Dbar = \iota|\Dbar$ reverses the orientation of $\Dbar$, $\tau$ acts
nontrivially on homology. Hence $[\tau] \in \pi_0(\SS_2)$ maps to the nontrivial
element of $\pi_0(\SS_1)$. $\tau^2$ can be deformed to the identity map by a
path $(h_s)_{0 \leq s \leq 1}$ in $\SS_1$, e.g.
\[
h_s(x,y) = 
\begin{cases}
(x,y) & |x+y| \leq \half,\\
\rho(2s(\pi + r(|x+y|))) & \text{otherwise.}
\end{cases}
\]
The derivative of $h_s$ along $\Delta$ is a rotation of $\nu\Delta$ with angle $2\pi s$.
It follows that $[\tau]^2 = \partial(r)$ where $r \in \pi_1(\GG^s)$ is a generator. In
view of \eqref{eq:exact-sequence} this shows that $\pi_0(\SS_2) \iso \Z$, with 
$[\tau]$ as a generator. \qed

In the terminology of \cite{seidel98b} $\tau$ is the
generalized Dehn twist along the Lagrangian sphere $\Dbar$. Since we can identify
$\mathrm{int}(T)$ with $S^2 \times S^2 \setminus \Delta$ in such a way that $\Dbar$
corresponds to the zero-section $Z \subset \cotangent$, it follows that
$\pi_0(\SS)$ is generated by the generalized Dehn twist along $Z$.

\proof[Proof of Corollary \ref{th:ss-dd}] Consider the diagram of induced maps
\[
\xymatrix{
1 \ar[r] &
\pi_1(\GG^s) \iso \Z \ar[r]^-{\partial} \ar[d]_{\alpha_1}&
\pi_0(\SS_2) \iso \Z \ar[r] \ar[d]_{\alpha_2} &
\pi_0(\SS_1) \iso \Z/2 \ar[d]_{\alpha_3} \\
\pi_1(\DD_1) \ar[r] &
\pi_1(\GG) \ar[r]^-{\partial'} &
\pi_0(\DD_2) \ar[r] &
\pi_0(\DD_1).
}
\]
$\alpha_3$ is injective, and since $\GG^s$ is a deformation retract of $\GG$,
$\alpha_1$ is an isomorphism. It remains to show that $\partial'  = 0$,
or equivalently, to find a loop $(\lambda_t)_{0 \leq t \leq 1}$ in $\DD_1$ whose
action on $\nu\Delta$ is a generator of $\pi_1(\GG)$. Such a loop is given by
$\lambda_t(x,y) = (R_y^{2\pi t}(x),y)$. \qed

\proof[Proof of Lemma \ref{th:ss-one}] A map $\psi \in \DD_1$ can act on
$H_2(S^2 \times S^2)$ in two ways: either trivially or, like $\iota$, by
interchanging the two generators. Let $\DD_+ \subset \DD_1$ be the
subgroup of those $\psi$ which act trivially, and $\SS_+$ the corresponding
subgroup of $\SS_1$. It is sufficient to
show that $\SS_+$ is contractible. To do this we use an argument due
to Gromov \cite[2.4.$A_1$]{gromov85} (see \cite{abreu97} for further
applications of this kind of argument). Let $\J$ be the space of those
$\o$-compatible almost complex structures on $S^2 \times S^2$
whose restriction to $T\Delta$ agrees with that of the standard complex
structure $J_0$. There is a map $A: \SS_+ \longrightarrow \J$
given by $A(\phi) = \phi^*(J_0)$. Gromov's argument provides a
map $B: \J \longrightarrow \SS_+$ such that $B \circ A = \id$;
since $\J$ is contractible, this proves that $\SS_+$ is also
contractible.

Recall that any $J \in \J$ determines, through its pseudo-holomorphic
curves, a pair $(\F_1,\F_2)$ of two-dimensional symplectic foliations of 
$S^2 \times S^2$, which are transverse to each other. The leaves of these 
foliations are two-spheres whose homology class is $A_1 = [S^2 \times \point]$ 
resp. $A_2 = [\point \times S^2]$. Since $\Delta$ is $J$-holomorphic and 
$[\Delta] \cdot A_i = 1$, any leaf of $\F_1$ or $\F_2$ intersects 
$\Delta$ transversely in a single point. It follows that there is a
unique $\phi \in \DD_+$ which maps $\F_1$ and $\F_2$ to the two 
standard foliations of $S^2 \times S^2$ by two-spheres. $\phi$ is
not symplectic, but since $\phi^*\o$ is tamed by $J$, one can use
Moser's technique to deform $\phi$ into a symplectic automorphism $\phi'$;
this can be done in such a way that $\phi' \in \SS_+$. Since $\phi =
\phi(J)$ depends smoothly on $J$, one can arrange that $\phi' = \phi'(J)$ 
also depends smoothly on $J$, and that $\phi'(J) = \phi(J)$ for all
those $J$ for which $\phi$ is already symplectic. Then one defines
$B$ by $B(J) = \phi'(J)$. \qed

\begin{remark} The higher-dimensional analogue of the problem which
we have considered presents some potentially more subtle features.
Let $T_n$ be the bundle of unit discs in $T^{*\!}S^n$,
and $\SS_n$ the group of symplectic automorphisms of $T_n$ which are
the identity near the boundary. Using a generalization of the argument
in \cite{seidel98b} it can be shown that $\pi_0(\SS_n)$ contains a 
class $[\tau]$ of infinite order. However, there is a possible source
for elements other than powers of $[\tau]$. The group $\pi_0(\Diff^+(S^n))$,
which is isomorphic to the group of exotic $(n+1)$-spheres for $n \geq 4$,
acts on $\pi_0(\SS_n)$ by conjugation, and it is not known whether this
action preserves $[\tau]$. \end{remark}

%------------------------------------------------------------------------
\providecommand{\bysame}{\leavevmode\hbox to3em{\hrulefill}\thinspace}

\end{document}